\documentstyle[amssymb,amsmath,eucal]{article}

\begin{document}


\def\K{{\Bbb K}}
\def\R {{\Bbb R }}
\def\C {{\Bbb C }}
\def\H {{\Bbb H }}

\def\U{{\rm U}}
\def\Sp{{\rm Sp}}
\def\const{{\rm const}}
\def\B{{\rm B}}
\def\O{{\rm O}}
\def\SO{{\rm SO}}
\def\SOS{{\rm SO}^*}
\def\GL{{\rm GL}}
\def\SL{{\rm SL}}
\def\SU{{\rm SU}}

\def\UU{{\cal U}}

\def\Gr{{\rm Gr}}
\def\Mat{{\rm Mat}_{p,q}}
\def\G{\Gr_{p,q}}

\def\ov{\overline}

\def\phi{\varphi}
\def\epsilon{\varepsilon}
\def\le{\leqslant}
\def\ge{\geqslant}

\def\tg{\tan}

\def\kvadrat{\hfill{$\boxtimes$}}

\newcounter{punct}
\newcounter{fact}

\def\punct{\addtocounter{punct}{1}\arabic{punct}. }
\def\fact{\addtocounter{fact}{1}\arabic{fact}.}

\begin{center}

{\large\bf  On Jordan angles and triangle inequality in Grassmannian}

\vspace{22pt}

Yurii A.Neretin%
\footnote{Supported by grants RFBR--98-01-00303,  NWO 047-008-009}

\end{center}

\begin{abstract}
Let $L,M,N$ be $p$-dimensional subspaces in $\R^n$.
Let $\phi_j$ be the angles between $L$ and $M$,
let $\psi_j$ be the angles between $M$ and $N$,
and $\theta_j$ be the angles between $L$ and $M$.
Consider the orbit of the vector $\psi\in\R^p$
with respect to permutations of coordinates and inversions 
of axises. Let $Z$ be the convex hull of this orbit.
Then $\theta\in\phi+Z$. We discuss similar theorems for other
symmetric spaces.

We also obtain formula for geodesic distance
on any invariant convex Finsler metrics
on classical symmetric space
\end{abstract}

\vspace{22pt}

We obtain a version of V.B.Lidskii  theorem
\cite{Lid1} on  spectrum of a sum of matrices for
arbitrary classical Riemannian semisimple symmetric space
\footnote{The result of the paper was announced in \cite{Ner2}}.

\smallskip

{\bf \punct  Grassmannians.} Fix  positive integers $p\le q$.
Consider the space $\R^{p+q}$ equipped with the standard
scalar product.
Denote by $\G$   the set of all $p$-dimensional linear subspaces in
$\R^{p+q}$. The orthogonal  group $\O(p+q)$ acts in
      $\R^{p+q}$  and hence it acts on $\G$. Obviously,
$$\G=\O(p+q)/\O(p)\times\O(q)$$

\smallskip

{\bf \punct Jordan angles.}
Let $L,M\in\G$. Consider an orthonormal basis $e_1,\dots, e_p\in L$
and an orthonormal basis $f_1,\dots,f_p\in M$. Consider the matrix
  $$\Lambda=\Lambda[L,M]$$
with the matrix elements $<e_i,f_j>$.
Denote by
$$\lambda_1\ge\dots \ge\lambda_p$$
the singular values%
\footnote{{\it The singular values} of a matrix $A$ are the eigenvalues
of the matrix $\sqrt{A^*A}$.}
  of the matrix $\Lambda$.
Obviously, the numbers
$$\lambda_j=\lambda_j[L,M]$$
don't depend on  choice of the bases
$e_1,\dots, e_p\in L$ and $f_1,\dots, f_p\in M$.

\smallskip

{\sc  Proposition \fact} {\it Let $L,M, L',M'\in\G$. The following
conditions are equivalent

\smallskip

{\rm i)} $\lambda_j[L,M]=\lambda_j[L',M']$    for all $j$

\smallskip

{\rm ii)} There exists an element $g\in\O(p+q)$
such that $gL=L', gM=M'$.   }

\smallskip

{\sc Proof.} The statement is obvious.

\smallskip

{\sc Proposition \fact} {\it  Consider $L,M\in\G$.
 There exist  orthonormal bases
$e_1,\dots, e_p\in L$ and $f_1,\dots, f_p\in M$
such that}
\begin{align*}
&<e_i,f_j>=0 \qquad\mbox{if}\qquad i\ne j \\
&<e_j, f_j>= \lambda_j
\end{align*}

{\sc Proof.} The statement is obvious.

\smallskip

{\sc Proposition \fact}
$$a)\qquad \lambda_k[L,M]=
          \max\limits_{\begin{array}{c}P\subset L\\ \phantom{P}\end{array}}
            \min\limits_{\begin{array}{c} v\in P\\ \|v\|=1\end{array}}
            \max\limits_{\begin{array}{c} w\in M\\ \|w\|=1\end{array}}
               <v,w>
$$
{\it where the first maximum is given over all $k$-dimensional subspaces
$P$ in $L$.}

$$b)\qquad \lambda_k[L,M]=
         \min\limits_{\begin{array}{c}Q\subset L\\ \phantom{P}\end{array}}
            \max\limits_{\begin{array}{c} v\in Q\\ \|v\|=1\end{array}}
            \max\limits_{\begin{array}{c} w\in M\\ \|w\|=1\end{array}}
               <v,w>
$$
{\it where the  minimum is given over all  subspaces  $Q$
having codimension $k$
in $L$.}

\smallskip

{\sc Proof.}   The statement is a corollary of the standard
minimax characterizations of eigenvalues and singular values, see
\cite{Lid2}, \cite{Bha}.

\smallskip

{\sc Proposition \fact}  {\it Denote by  $\Pi_M$ the orthogonal
projector to the subspace  $M$. Then the numbers
$\lambda_j[L,M]$ are the singular values of the operator
$$\Pi_M: L\to M$$.}

\smallskip

{\sc Proof.} The statement is obvious.

\smallskip

{\sc Lemma \fact} {\it
Let $u_1, u_2,\dots,u_p$ be a
{\rm(} nonorthogonal {\rm)} basis in L.
Let $v_1, v_2,\dots,v_p$ be a {\rm(} nonorthogonal {\rm)}
 basis in M.
Denote by $U$ the matrix with the matrix elements $<u_i,u_j>$, denote by
$V$ the matrix with the matrix elements $<v_i,v_j>$,
denote by $W$ the matrix with the matrix elements $<u_i,v_j>$.
Then the numbers $\lambda_j^2[L,M]$ coincides with
eigenvalues of the matrix}
\begin{equation}
 U^{-1}W V^{-1} W^t
\end{equation}

{\sc Proof.} The statement is obvious.

{\it The angles} or {\it stationary angles}%
\footnote{Other terms are {\it complex distance}
or {\it compound distance}. Models of symmetric spaces
given in \cite{Ner1} and \cite{Ner4} shows that these invariants
for all classical Riemannian symmetric spaces are really angles.}
 (K.Jordan, 1875) $\psi_1\le\psi_2\le\dots\le \psi_p$
 between the subspaces
$L, M\in\Gr_{p,q}$ are defined by
$$\psi_1=\Psi_1[L,M]:=\arccos\lambda_1, \,
 \dots,\, \psi_p=\Psi_p[L,M]: =\arccos\lambda_p $$
Obviously, $0\le \psi_j\le \pi/2$. We also will use the notation
$$\Psi[L,M]=(\Psi_1[L,M],\dots,\Psi_p[L,M])$$

{\sc Remark.} For all $j$ we have $\Psi_j[L,M]=\Psi_j[M,L]$.

\smallskip

{\bf \punct The result of the paper.}
Denote by $W_p$ (Weyl group)
 the group of all transformations of $\R^p$
generated by permutations of the coordinates and  by the
transformations
$$(t_1,\dots,t_p)\mapsto (\sigma_1 t_1,\dots, \sigma_p t_p)
 \qquad\mbox{where}\qquad\sigma_j=\pm 1
$$

\smallskip

{\sc Theorem A.} {\it Let $\ell(x)$ be a $W_p$-invariant
 norm in $\R^p$. Then the function
$$
d(L,M):=\ell(\Psi_1[L,M], \dots, \Psi_p[L,M])
$$
is  an $\O(p+q)$-invariant metric on $\G$.}

\smallskip

{\sc Remark.} The geodesic distance in $\Gr_{p,q}$ associated with
the $\O(p)\times\O(q)$-invariant Riemannian metrics is given by the formula
$${\rm dist}(L,M)=\sqrt{\Psi_1[L,M]^2+ \dots+ \Psi_p[L,M]^2}$$

{\sc Theorem B.} {\it Let $L,M,N\in\G$. Let $\phi_j=\Psi_j[L,M]$,
$\psi_j=\Psi_j[M,N]$, $\theta_j=\Psi_j[L,N]$ be the angles.
Denote by ${\cal Z}$ the convex hull of the $W_p$-orbit of
the vector $(\psi_1,\dots,\psi_p)\in \R^p$. Denote by
$\cal U$ the shift of $\cal Z$ by the vector
$(\phi_1,\dots,\phi_p)$.

Then there exists a vector $(\theta_1^\circ,\dots,\theta_p^\circ)\in \cal U$
such that the collection of numbers
$(\cos \theta_1^\circ,\dots,\cos\theta_p^\circ)$  coincides
up to permutation with the collection of numbers
$(\cos \theta_1,\dots,\cos\theta_p)$.}

\smallskip

{\bf \punct Infinitesimal angular structure.}
For $L\in\G$ we denote by $T_L(\G)$ the tangent space to $\G$
at the point $L$. It is natural to identify
elements $\xi\in T_L(\G)$ with operators $H$ from $L$ to
the orthogonal complement $L^\bot$. We denote
by
$$\rho_1[L;H]\le\dots\le \rho_p[L;H]$$
the singular values of the operator $H:L\to L^\bot$.

\smallskip

{\sc Lemma \fact } {\it Let $L,L'\in\G$, $H\in T_L(\G)$ , $H'\in T_{L'}(\G)$.
The following conditions are equivalent }

\smallskip

i) {\it $\rho_j[L;H]= \rho_j[L';H']$ for all $j$}

\smallskip

ii) {\it There exists an operator $g\in\O(p+q)$ such that
$gL=L'$, $gH=H'$.}

\smallskip

{\sc Proof.} The statement is obvious.

\smallskip

{\sc Remark.} The $\O(p+q)$-invariant Riemannian metric in $\G$
is
$${\rm tr}\,\, H^tH=\sum \rho_j^2[L;H]$$

\smallskip

{\bf \punct Relations between angles and infinitesimal
angular structure.}  The following statement is obvious.

{\sc Proposition \fact} {\it Let $M(t)$ be a smooth $C^\infty$-curve
in $\G$. Then
$$\lim_{\epsilon\to +0}\frac{\Psi_j[M(a+\epsilon),M(a)]}{\epsilon}
         =\rho_j[M(a);M'(a)]$$
where $\lim\limits_{\epsilon\to +0}$ denotes the right limit at $0$.}

\smallskip

{\bf \punct Infinitesimal variation of angles.}
 Let $L,M\in \G$. Let $e_j\in L$,
$f_j\in M$ be orthonormal bases satisfying the conditions
\begin{align*}
&<e_i,f_j>=0 \qquad\mbox{if}\qquad i\ne j \\
&<e_j, f_j>= \cos \psi_j
\end{align*}
{\it Assume $\psi_j$ be pairwise different.}

Consider an orthonormal basis $r_1,\dots,r_q\in M^\bot$ such that
for all $j\le p$ the vectors
$e_j$, $f_j$, $r_j$ span 2-dimensional plane and $f_j$ is situated
in the angle between $e_j$ and $r_j$.

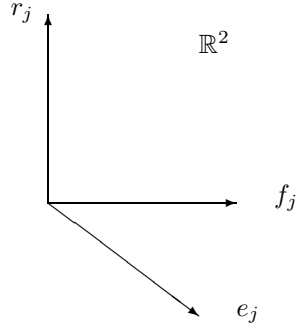
\unitlength5mm
\begin{picture}(7,9)
\put(5,3){
\put(0,0){\vector(0,1){5}}
\put(0,0){\vector(1,0){5}}
\put(0,0){\vector(4,-3){4}}
\put(6,0){$f_j$}
\put(5,-3){$e_j$}
\put(-1,5){$r_j$}
\put(4,4){$\R^2$}
}
\end{picture}

We have
$$e_j=f_j\cos\psi_j-r_j\sin\psi_j$$

Let $H:M\to M^\bot$ be a tangent vector to $\G$ at the point $M$.
Let $h_{ij}$ be the matrix elements of $H$ in the bases
$f_1,\dots,f_p$ and $r_1,\dots,r_q$.

Consider a $C^\infty$-smooth curve $M(\epsilon)\in\G$ such that
$$M(0)=H; \qquad M'(0)=H$$

\smallskip

{\sc Proposition \fact} {\it Assume $\psi_0\ne 0$, $\psi_p \ne \pi/2$,
and $\psi_{j+1}\ne \psi_{j}$ for all $j$. Then                  }
\begin{equation}\frac{d}{d\epsilon}\Psi_j[L,M(\epsilon)]\Bigr|_{\epsilon=0}
                 = h_{jj}
\end{equation}


{\sc Proof.} Denote by $f_j(\epsilon)$ the unique
vector in $M(\epsilon)$ having the form
$$f_j(\epsilon)=f_j+\sum a_{jk}(\epsilon)r_k$$
Obviously,
$$a_{jk}(\epsilon)=\epsilon h_{jk}+O(\epsilon^2)$$
For all $i,j\le p$ we have
$$<f_i(\epsilon),f_j(\epsilon)>=\left\{
   \begin{array}{c} O(\epsilon^2),\qquad \mbox{if}\qquad i\ne j \\
                    1+O(\epsilon^2),\qquad \mbox{if}\qquad i= j
   \end{array}\right.
$$
and
$$<e_i,f_j(\epsilon)>=\left\{
   \begin{array}{c}
    -\epsilon h_{ij}\sin \psi_i+ O(\epsilon^2),\qquad \mbox{if}\qquad i\ne j \\
                    \cos(\psi_j+\epsilon h_{jj})+O(\epsilon^2),\qquad \mbox{if}\qquad i= j
   \end{array}\right.
$$
Now we are ready to wright matrix (1) for the subspaces $L$, $M(\epsilon)$
$$
\begin{pmatrix}
\cos^2 (\psi_1+\epsilon h_{11})+O(\epsilon^2)
                                   & O(\epsilon)&\dots & O(\epsilon)\\
O(\epsilon)& \cos^2 (\psi_2+\epsilon h_{22})+O(\epsilon^2) &\dots & O(\epsilon) \\
\vdots & \vdots & \ddots & \vdots \\
O(\epsilon)  & O(\epsilon) &\dots &
            \cos^2 (\psi_p+\epsilon h_{pp})+O(\epsilon^2)
\end{pmatrix}
$$
This implies required statement.

\smallskip

{\sc Remark.} Let us fix $L\in\G$.

  a) The set of all $M\in\G$ such that $\Psi_0[L,M]=0$ has codimension
$q-p+1$.

  b) The set of all $M\in\G$ such that $\Psi_p[L,M]=\pi/2$ has codimension
1.

  c) The set of all $M\in\G$ such that $\Psi_{j+1}[L,M]=\Psi_{j}[L,M]$ has
 codimension 2.

\smallskip

{\bf \punct Preliminaries.} A $p\times p$ matrix $A$ is called
{\it bistochastic} if for all $k$, $l$
$$\sum_i a_{ik}=1;\qquad \sum_j a_{lj}=1$$

We  say that a real $p\times p$ matrix $A$ is {\it quasistochastic}
if for all $k$, $l$
$$\sum_i |a_{ik}|\le 1;\qquad \sum_j |a_{lj}|\le 1$$

\smallskip

{\sc Proposition \fact} (Birkhoff) {\it The set of all bistochastic matrices
is the convex hull of matrices of permutations}%
\footnote{i.e. matrices consisting of 0 and 1 and having strictly
one 1 in each column and each row.}

\smallskip

See \cite{Lid2}, \cite{Bha}.

\smallskip

{\sc Lemma \fact} {\it  The set of all quasistochastic matrices
is the convex hull of the group $W_p$.}

\smallskip

{\sc Proof.} It is sufficient to describe extremal points of
the set of all quasistochastic matrices.

 a) Obviously, for any extremal point $A$
\begin{equation}
\sum_i |a_{ik}|= 1;\qquad \sum_j |a_{lj}|= 1
\end{equation}

 b) Let a matrix $A$ satisfies condition (3). Assume
the matrix $|a_{ij}|$  be not an extremal point of the set
 of stochastic matrices.  Then $A$ is not an extremal point of the
set of quasistochastic matrices.

  Hence any extremal point of the set of quasistochastic matrices
is an element of $W_p$.

\smallskip

{\sc Example.}
Let $U=(u_{ij})$ and $V=(v_{ij})$ be matrices with norm $\le 1$.
 Then
the matrix $W=(u_{ij}v_{ij})$ is quasistochastic.

          \smallskip

{\sc Lemma \fact}%
\footnote{This is a minor variation of Fan Ky theorem, see
\cite{Lid2}, \cite{Bha}}.
{\it Let $A$ be a real $p\times q$ matrix. Let
$\lambda=(\lambda_1,\dots,\lambda_p)$ be its singular
values. Then the convex hull of the $W_p$-orbit of $\lambda$
contains the vector $(a_{11},\dots, a_{pp})$.}

\smallskip

{\sc Proof.} Indeed, the matrix $A$ can be represented in the form
$$A=U \begin{pmatrix}\lambda_1 &0 &\dots&0&0&\dots&0\\
                     0&\lambda_2 &\dots&0&0&\dots&0 \\
                   \vdots  &\vdots&\ddots     &\vdots&\vdots&\vdots&\vdots \\
                    0&0 &\dots& \lambda_p &0&\dots&0
       \end{pmatrix} V^t;\qquad U\in\O(p),\,\, V\in \O(q)
$$
Hence
$$
\begin{pmatrix}a_{11}\\a_{22}\\ \dots \\ a_{pp}\end{pmatrix} =
        \begin{pmatrix} u_{11}v_{11} &  u_{12}v_{12}&\dots& u_{1p}v_{1p} \\
                        u_{21}v_{21} &  u_{22}v_{22} &\dots& u_{2p}v_{2p}  \\
                        \dots & \dots &\dots&\dots \\
                        u_{p1}v_{p1} &  u_{p2}v_{p2} &\dots& u_{pp}v_{pp}
                                                       \end{pmatrix}
         \begin{pmatrix} \lambda_1 \\ \lambda_2 \\ \dots \\ \lambda_p
         \end{pmatrix}
$$
Then we apply Lemma 10.

{\bf \punct Proof of Theorems A--B.}
Fix nonnegative numbers
$$a_1\le\dots\le a_p$$
Fix arbitrary orthonormal basis
\begin{equation}
e_1,e_2\dots,e_p,
f_1,f_2\dots,f_p, r_1,\dots, r_{q-p}\in \R^{p+q}
\end{equation}
Consider the subspace $L_a(s)\in \G$ spanned by the vectors
$v_1(s),\dots,v_p(s)$ given by the formula
$$v_j(s)=\cos (a_js)e_j+\sin(a_js)f_j$$
We obtain a curve $L_a(s)$ in $\G$.

We say  that a curve  $\gamma(t)$ in $\G$
is a {\it $H$-curve} if in some orthonormal basis it
has the form  $L_a(s)$.  We say that the numbers $a_j$ are
the {\it invariants}
of the $H$-curve  $\gamma$.

We say that points $L_a(s)$, $L_a(t)$ on a $H$-curve are
{\it sufficiently near} if
$$a_p|s-t|\le \pi/2$$

The following statements are obvious

\smallskip

{\sc Lemma \fact} {\it
Consider sufficiently near points $L(s_1)$, $L(s_2)$, $L(s_3)$
on $H$-curve. Assume $s_1<s_2< s_3$.
 Then for all $j$}
$$\Psi_j [L(s_1),L(s_2)]+ \Psi_j [L(s_2),L(s_3)]
      = \Psi_j [L(s_1),L(s_3)]
$$

\smallskip

{\sc Lemma \fact} {\it  Let $L, M\in \Gr_{p,q}$ and
$\Psi_p[L,M]<\pi/2$. Then there exists
the unique $H$-curve $\gamma(s)$  joining $L$, $M$ such that
$L$, $M$ are sufficiently near points on $\gamma(s)$.}

\smallskip

Consider points $L,M,N\in \G$ having a general position.
Denote by $\theta_j$ the angles between $M$ and $N$.
Consider the $H$-curve $\gamma(t)$ such that $\gamma(0)=M$,
$\gamma(1)=N$ and $M$, $N$ are sufficiently near
points of the curve $\gamma(t)$. Then the invariants of the $H$-curve
$\gamma(s)$ are $\theta_1,\dots,\theta_p$.
 We assume, that for each $s\in [0,1]$
\begin{gather}
\Psi_0[L,M(s)]\ne 0,\qquad \Psi_p[L,M(s)]\ne \pi/2 \\
\Psi_j[L,M(s)] \ne \Psi_{j+1}[L,M(s)] \qquad \mbox{for all}\quad j
\end{gather}

 Denote by $\cal Z$ the convex hull of $W_p$-orbit
of the vector $\theta$.

 By Proposition 8 and Lemma 11,
we have
$$
\Psi[L,\gamma(1/n)]\in \Psi[L,M]+\bigl(\frac{1}{n} +O(\frac1{n^2})\bigr)
{ \cal Z},
\qquad n\to\infty
$$
In the same way,
\begin{align*}
&\Psi[L,\gamma(2/n)]
       \in \Psi[L,\gamma(1/n)]+\bigl(\frac{1}{n} +O(\frac1{n^2})\bigr) \cal Z  \\
&\Psi[L,\gamma(3/n)]
       \in \Psi[L,\gamma(2/n)]+\bigl(\frac{1}{n} +O(\frac1{n^2})\bigr) \cal Z\\
&.\qquad .\qquad .\qquad .\qquad .\qquad .\qquad .\qquad.
\end{align*}
and $O(1/n^2)$ are uniform in  $k/n$.
Hence
$$
\Psi[L,\gamma(t)]\in
\Psi[L,M] + t\bigl(1+O(\frac1n)\bigr) {\cal Z}; \qquad n\to\infty
$$
and hence
\begin{equation}
          \Psi[L,\gamma(t)]\in \Psi[L,M] + t {\cal Z};
\end{equation}

Assume that there exists the unique value $\widetilde s$ that
doesn't satisfies the conditions (5)--(6). Consider a small $\delta$. Then
for $t>\widetilde s$ we have

\begin{align*}
&\Psi[L,\gamma(\widetilde s-\delta)]\in \Psi[L,M] +
              (\widetilde s-\delta) {\cal Z}\\
&\Psi[L,\gamma(\widetilde s-\delta)] \qquad\mbox{is close to}\qquad
   \Psi[L,\gamma(\widetilde s+\delta)] \\
&\Psi[L,\gamma(t)]\in
        \Psi[L,\gamma(\widetilde s+\delta)] +
             (t- \widetilde s+\delta) {\cal Z}
\end{align*}
and we again obtain (7).

A $H$-curve of general position contains a finite number of points
$\widetilde s$ that don't satisfy conditions (5)--(6).
Hence we can repeat our arguments.

This proves Theorem B.

Consider a $W_p$-invariant norm $\ell(\cdot)$ on $\R^p$.
By Lemma 11 for any $x\in\cal Z$
$$\ell(x)\le\ell(\theta)$$
and this finishes the proof of  Theorem A.

\smallskip

{\sc  Corollary \fact} {\it $H$-curves are geodesics in any
metrics having the form
\begin{equation}
d(L,M)=\ell(\Psi_1[L,M],\dots,\Psi_p[L,M])
\end{equation}}

\smallskip

{\sc Remark.} If the sphere $\ell(x)=1$ in $\R^n$ doesn't contain a segment,
then this geodesics is unique.

\smallskip

{\sc  Corollary \fact}  {\it Consider the Finsler metric on $\Gr_{p,q}$
given by the formula
$$F(L,H)=\ell(\rho_1[L,H],\dots, \rho_1[L,H]);\qquad L\in\Gr_{p,q},H\in T_L$$
Then the associated geodesic distance is given by the expression \rm (8).}

{\bf \punct Other symmetric spaces.}
In  \cite{Ner4} it was explained that arbitrary classical
compact symmetric space is a Grassmannian in real, complex or quaternionic
linear space.                 This
allows to translate literally our results to
all classical  compact Riemannian
symmetric spaces
\begin{align}
&\U(n) \times \U(n) /\U(n);\quad \U(n)/\O(n); \qquad \U(2n)/\Sp(n)\\
&\U(p+q)/\U(p)\times\U(q);\quad
     \O(2n)/\U(n);\quad \O(p+q)/\O(p)\times\O(q);
                             \quad \O(n)\times \O(n)/\O(n);\notag\\
&\Sp(p+q)/\Sp(p)\times\Sp(q);\quad \Sp(n)/\U(n);\quad \Sp(n)\times \Sp(n) /\Sp(n)
   \notag
\end{align}
Three series of the type $A$ (i.e (9))
 slightly differs from others:
we have to replace the group $W_p$ by the symmetric group.

In the same way, all classical Riemannian noncompact symmetric spaces
are open domains in Grassmannians
(see \cite{Ner1}). This
 allows
to extend our results  to all classical
 Riemannian noncompact symmetric spaces
\begin{align*}
&\GL(n,\C)/\U(n);\quad \GL(n,\R)/\U(n);\quad \GL(n,\H)/\U(n)\\
&\U(p,q)/\U(p)\times\U(q);\quad
     \SO^*(2n)/\U(n);\quad \O(p,q)/\O(p)\times\O(q);
                             \quad \O(n,\C)/\O(n);\\
&\Sp(p,q)/\Sp(p)\times\Sp(q);\quad \Sp(2n,\R)/\U(n);\quad \Sp(n,\C) /\Sp(n)
\end{align*}

{\bf \punct Some examples.} a) {\it The space $\GL(n,\C)/\U(n)$.}
We realize points of the space as positive definite $n\times n$ complex
matrices.  The group $\GL(n,\C)$ acts on this space by the transformations
$$L\mapsto gLg^*,\qquad g\in \GL(n,\C)$$

 The {\it angles}\footnote{hyperbolic angles}
 $\Psi_j[L,M]$
 between points $L$ and $M$  are the
solutions of the equation
$$\det(L-e^\psi M)=0$$
Denote by $\Psi[L,M]$ the vector  $(\Psi_1[L,M],\dots,\Psi_n[L,M])$.

Let $L$, $M$, $N$ be points of our space.
Consider all vectors in $\R^n$ that can be obtained from
$\Psi[L,M]$   by permutations of coordinates.
 Denote by  $\cal Z$  their convex hull.
Then
\begin{equation}
\Psi[L,N]\in \Psi[L,M]+\cal Z
\end{equation}

b) {\it Original Lidskii theorem.} Consider the space $\cal S$ of hermitian
$n\times n$ matrices. This space also is a (nonsemisimple) symmetric space.
The group of isometries is the group of transformations
$$X\mapsto UXU^*+ A \qquad\mbox{where}\qquad U\in\U(n), A\in \cal S$$
Let $X,Y\in\cal S$. The analogy of angles are the eigenvalues of $X-Y$.
The analogy of Theorem B is the original Lidskii theorem \cite{Lid1}.

The space $\cal S$ can be identified with the tangent space to
$\GL(n,\C)$ at the point $1$. For $X,Y\in \cal S$ we define matrices
$$L=1+\epsilon A,\qquad M=1+\epsilon B\quad \in \GL(n,\C)/\U(n)$$
where $\epsilon$ is small. Then  the angles between  $L$ and $M$
have the form $\epsilon \lambda_j$, where $\lambda_j$ are the
eigenvalues of $X-Y$. Hence the inclusion (10) implies Lidskii
theorem.

 Lidskii theorem on singular values
of sum of two matrices  corresponds to the triangle inequality
in a tangent space to $\U(p+q)/\U(p)\times\U(q)$ or $\U(p,q)/\U(p)\times\U(q)$

\smallskip

c) {\it The space $\Sp(2n,\R)/\U(n)$.} This spaces can be realized as
the space of symmetric $n\times n$ complex matrices with norm $<1$
(see for instance \cite{Ner3},5.1,6.3).
For two points $T,S$ we define the
 expression
\begin{equation}
\Lambda[T,S]=(1-TT^*)^{-1/2}(1-TS)(1-SS^*)^{-1/2}
\end{equation}
Let $\lambda_j$ be its singular values. Then the
{\it hyperbolic angles} between $T$ and $S$ are given
by the formula
\begin{equation}
\psi_j={\rm arcosh}\,\, \lambda_j
\end{equation}
The analogy of Theorem B is given by the formula (10).

\smallskip

d) {\it Arazy norms.}
Denote by $V_{fin}$  the space
of finite real  sequences $x=(x_1,\dots,x_N,0,0,\dots)$.
Consider a norm $\ell$ on a space $V_{fin}$
 satisfying the conditions

  -- $\ell$ is invariant with respect to permutations of coordinates

  -- $\ell$ is invariant with respect to the transformations
     $(x_1,x_2,\dots)\mapsto (\sigma_1 x_1,\sigma_2 x_2,\dots)$, where
$\sigma_j=\pm 1$

-- if $x^{(j)}$ converges to $x$ coordinate-wise and
$\ell(x^{(j)})$ converges to $\ell(x)$, then $\ell(x^{(j)}-x)$ converges
to 0 (an equivalent formulation: the $\ell$-convergence on the sphere
$\ell(x)=1$ is equivalent to  the coordinate-wise convergence).

 Let $V_\ell$ is the completion of $V_{fin}$ with respect to the
norm $\ell$.

A compact operator $A$ in a Hilbert space is an element of {\it  Arazy class}
(see \cite{Ara}) $C_\ell$ if the sequence of its singular values is an element of $V_\ell$.
Consider a space  ${\rm B}_\ell$  (operator ball)
 of all compact operators $T$ in the Hilbert space satisfying the conditions

-- $T \in C_\ell$

-- $\|T\|<1$, where $T$ denotes the standard norm of a operator
in a hilbrt space

-- $T=T^t$

We define the angles $\Psi_j(T,S)$ in  ${\rm B}_\ell$
  by  formulas (11)--(12).
We define the distance in  ${\rm B}_\ell$  by
$$d_\ell(T,S)=\ell (\Psi_1(T,S),\Psi_2(T,S),\dots)$$

{\sc Proposition \fact} a) {\it $d_\ell(T,S)$ is  a metric.}

b) {\it The space   ${\rm B}_\ell$      is complete with respect
to the metric $d_\ell(T,S)$.}

\smallskip

The statement a) can be easyly obtained from Theorem A by a limit
considerations. For a proof of the statement B   see
\cite{Ner3}, 8.6.3.

\smallskip

{\bf \punct Some references.}
a) For matrix inequalities  see for instance \cite{Lid2}, \cite{Bha}.

\smallskip

b) The formula for the  distance in a symmetric space
 associated with the invariant
Riemannian metrics was obtained in \cite{Kli}

\smallskip

c) Our Theorem B for the unitary group $\U(n)= \U(n)\times \U(n)/\U(n)$
is Nudelman--Shvartsman theorem \cite{NS}

\smallskip

d) Generalization of Fan Ky theorem to arbitrary
simple Lie algebras was obtained in \cite{Kos}.

e) Let $G$ be a simple Lie group, $K$ be its maximal
compact subgroup and $K\setminus G/K$
be the hypergroup of $K$-biinvariant subsets in $G$
with convolution product. The problem about triangle
inequality in Grassmannians is related to a classical
problem on  structure of the hypergroup
$K\setminus G/K$, see
for instance \cite{FK}, \cite{Rou},
\cite{DRW}.

\smallskip

f) Some nonstandard geometries on groups are discussed
in \cite{BP}.

\smallskip

g) Some applications of geometry of angles
 are contained in \cite{Ner3}, 6.3.

\smallskip

h) {\sc Conjecture.} I think that complete triangle inequality
for angles
coincides with Horn--Klyachko inequalities
\cite{Kly}, see also \cite{LidB} and some comments in
\cite{Zel}.

{\sc  Independent University of Moscow,
Bolshoj Vlas'evskij per., 11,
Moscow, 121002}

neretin@main.mccme.rssi.ru

\end{document}